\documentclass[12 pt]{amsart}

\usepackage{amssymb,amsfonts,graphicx}
\usepackage[all,arc]{xy}
\usepackage{enumerate}
\usepackage{mathrsfs}
\usepackage{enumitem} 
\newtheorem{theorem}{Theorem}[section]

\newtheorem{lemma}[theorem]{Lemma}

\theoremstyle{definition}

\newtheorem{notation}[theorem]{Notation}

\theoremstyle{remark}
\newtheorem{remark}[theorem]{Remark}


\numberwithin{equation}{section}

\title[FAMILIES OF SUBHARMONIC FUNCTIONS]{FAMILIES OF SUBHARMONIC FUNCTIONS AND SEPARATELY SUBHARMONIC FUNCTIONS}

\author{MANSOUR KALANTAR}

 \email{kalantarm@smccd.edu} 
\date{}

\begin{document}

\begin{abstract}
We prove that the upper envelope of a family of subharmonic functions defined on an open subset of $\mathbb{R}^{N}$, $(N\geq2)$, that is finite every where, is locally bounded above outside a closed nowhere dense set with no bounded components. Then we conclude as a consequence that a separately subharmonic function is subharmonic outside a closed nowhere dense set  with no bounded components. It generalizes a result due to Cegrell and Sadullaev.

\end{abstract}

\maketitle

\section{INTRODUCTION}

Let $u$ be the upper envelope  of a family   of subharmonic functions defined on an open set of $ \mathbb{R}^{N} $ ( $N\geq2$). A famous  result by Cartan \cite[Theorem 5.7.1.]{Armitage Gardiner book} stats that if $u$ is locally bounded above, then the upper regularization  $u^{*}$ of $ u $ is subharmonic and equals $u$ \textit{quasi everywhere}, i.e., outside a polar set.  In particular, $u=u^{*}$ \textit{almost everywhere}. It seems natural to try to characterize the  singular set of $u$, i.e., the set of points where $u$ fails to be locally bounded above. 

If the  functions   are  continuous and $u$ is real valued everywhere, it is an easy consequence of Bair spaces that $u$'s singular set is  a closed nowhere dense set \cite[Theorem 2. pg 194]{Bourbaki}. In the first part of this paper we generalize this result by relaxing the continuity condition  (Theorem 3.1); we use  the concept of thinness  and some results from the Perron-Wiener-Brelot (PWB) approach to the Dirichlet problem. Then, we prove under some natural condition that  conversely such a set  is the singular set of the supremum of a family of subharmonic functions (Theorem 3.5).

The second part of this paper is  an application of this (and related)  results to separately subharmonic functions; i.e., a function  $u(x,y)$  defined on an open set of $\mathbb{R}^{p}\times\mathbb{R}^{q} $ such that $u(a,.)$ and $u(,b)$ are subharmonic  for all fixed $a$ and $b$.

It is well-known that separately subharmonic functions are not necessarily (jointly) subharmonic, thanks to a  counter example by  J. Wiegerinck (\cite{Wiegerinck}, see also \cite{Wiegerinck Zeinstra}). However, by the work of V. Avanissian \cite{Avanissian} ( see also P. Lelong \cite{Lelong 1}), we know that these functions are subharmonic if and only if they are locally bounded above. Different authors replaced the boundedness in this result by weaker conditions and guaranteed the same conclusion: Arsove \cite{Arsove} replaced boundedeness by having a locally integrable majorant, and J. Riihentaus by having a locally $L^{p}$ majorant.   D. Armitage and S. Gardiner \cite{Armitage Gardiner} give an "almost charp" condition that includes all of the previous results. In \cite{Cegrell Sadullaev},  U. Cegrell and A. Sadullaev prove that such functions are subharmonic outside a product of closed nowhere dense set (singular set), if moreover they are harmonic with respect to one of the variables. 

We generalize Cegrell and Sadullaev's  result by relaxing the condition of harmonicity with respect to one of the variables (Theorem 4.1 below). We also show that the singular set of a separately subharmonic function has no bounded components, and in particular, no isolated points. At the end we prove the converse of this result under a mild extra condition (Theorem 4.2.).

\subsection*{Acknowledgments}  We wish to express out gratitude to Professor Ahmed Zeriahi who introduced the author to the fascinating world of potential theory, quiet some times ago.  It is a pleasure to thank Professor Azimbay Sadullaev for his help and advice during the past years to present. This is also a pleasure to thank Professor Juhani Riihentaus for his generous help to the author during last long years.

\section{DEFINITIONS AND PRELIMINARIES} 

\begin{notation}
 For a set $E\subset\mathbb{R}^{N}$, the closure, interior, and boundary of $E$ are noted $\overline{E}$, int$(E)$, and $\partial E$, respectively. The ball of center $ x $ and radius $ r>0 $ will be noted $ B(x,r) $. The Poisson kernel is noted $ K(\zeta,.) $ and $ \sigma $ designates the $ (N-1)$-dimensional Lebesgue measure on the boundary of  balls. The distance of a point $x$ and a set $E$ is noted $\text{d}(x,E)$; i.e.,
 $$\inf\lbrace \mid x-\zeta\mid: \zeta\in E\rbrace.$$
 \end{notation}

 For reader's convenience we summarize below some results from the classical potential theory that we will be using and can be found in \cite{Armitage Gardiner book}, \cite{Brelot} and \cite{Hayman Kennedy}. In what follows  $\Omega$ will be a bounded open subset of $\mathbb{R}^{N}$ for $N\geq2$. 
 
 We recall that an upper semi-continuous function $u:\Omega\rightarrow [-\infty,+\infty)$ is called subharmonic, if $u\not\equiv-\infty$ and  for all ball $B(x,\rho)$ relatively compact in $\Omega$, 
$$u(x)\leq \frac{1}{\sigma_{N}\rho^{N-1}}\int_{\partial B(x,\rho)}u(\zeta)d\sigma.$$

Let  $f$ be a continuous function on $\partial\Omega$. The harmonic extension of $f$ from $ \partial \Omega $ to $\Omega$ is the harmonic function defined by
\begin{equation}
H_{f}(x)=\int_{\partial \Omega}f(\xi)d\mu_{x}^{\Omega}(\xi),\label{2.1}
\end{equation}
where $\mu_{x}^{\Omega}$ is the harmonic measure at $x\in \Omega$. A point $ \zeta\in\partial\Omega $ is said to be a  regular (for Dirichlet problem) boundary point of $ \Omega $ if 
\begin{equation}
\lim_{\substack{x\rightarrow\zeta\\(x\in\Omega)}}H_{f}(x)=f(\zeta) \label{3018}
\end{equation}
for all function $ f $ continuous on $\partial \Omega $. Points for which the above equality does not hold are called irregular (for Dirichlet problem) boundary points of $ \Omega $.

An open set $\omega$ in $\mathbb{R}^{N}$ is said to be Greenian, if for each $y$ in $\omega$, the function defined by
           \begin{displaymath}
U_{y}(x) = \left\{
\begin{array}{lr}
-\log|x-y| & (x\not=y;N=2)\\
|x-y|^{2-N} & (x\not=y;N\geq3)\\
+\infty &  (x=y)
\end{array}
\right.
\end{displaymath}
has a subharmonic minorant on $\omega$. If $N\geq3$, all open sets of $\mathbb{R}^{N}$ are Greenian. For $\omega\subset \mathbb{R}^{2}$, if $\mathbb{R}^{2}\setminus\partial \omega$ is not connected, then $\omega$ is Greenian; in particular, all bounded open sets in  $\mathbb{R}^{2}$ are Greenian. See \cite[Definition 4.1.1. and Theorem 4.1.2]{Armitage Gardiner book}.

\begin{theorem}\label{1001}
 Let $ \omega $ be  Greenian. The set of irregular points (for the Dirichlet problem) of $\partial\omega$ is polar.

\end{theorem}
See \cite[Theorem 6.6.8]{Armitage Gardiner book}.

A set $ E\subset\partial \Omega $ is called negligible for $ \Omega $, if 
$$ d\mu_{x}^{\Omega}(E)=0 $$
for all $ x\in \Omega $. 
Negligible sets can by characterized by the notion of thinness. A set $ E $ is said to be thin at a point $ \zeta $ if $ \zeta $ is not a fine (with respect to the fine topology) limit point of $ E. $ The following two theorems, used in the proof of Theorem 3.1, give the relation between thinness and negligibility.

 \begin{theorem}\label{1002}
 Let $ \zeta $ be a limit point of a set $ E $. The set $ E $ is thin at $ \zeta $ if and only if there exists a subharmonic function $ v $ on a neighborhood of $ \zeta $ such that 
 $$ \limsup_{\substack{x\rightarrow\zeta\\(x\in E)}}v(x)<v(\zeta). $$
\end{theorem}
See \cite[ Theorem 7.2.3]{Armitage Gardiner book}.

\begin{theorem}\label{1000}
 Let $ \omega $ be Greenian.
\begin{itemize}
\item[(i)]
 The set $ \lbrace \zeta\in \partial \omega: \omega \text{ is thin at }\zeta \rbrace$ is negligible for $ \omega $.
 \item[(ii)]If $ E $ is a relatively open subset of $ \partial\omega $ which is negligible, then each point of $ E $ is irregular and the set $ E $ is polar.
 \end{itemize}
\end{theorem}
 
 See  \cite[Theorem 7.5.4]{Armitage Gardiner book} for part (i) and   \cite[Theorem 6.6.9-(i)]{Armitage Gardiner book} for part (ii).

 \begin{theorem}\label{th2.7} 
 Let $F$ be a compact subset of an open set $\Omega$ in  $ \mathbb{R}^{N}$ such that every bounded component of $\mathbb{R}^{N} \setminus F$ contains a point of $\mathbb{R}^{N}\setminus\Omega$. If $h$ is harmonic on an open set containing $F$ and if $\varepsilon>0$, then there exists a function $H$ harmonic on $\Omega$ such that $|h-H|<\varepsilon$ on $F$.
 \end{theorem} 
 See \cite[Theorem 2.6.4.]{Armitage Gardiner book}.
 
\vspace*{1cm}

In the proof of Theorem 3.1, we use the following result, sometimes called the Baire theorem: \textit{ The necessary and sufficient condition that a topological space $ X $ be a Baire space is that if $ X=\bigcup_{k}F_{k} $ and each $ F_{k} $ is closed, then $  \bigcup_{k}\mathrm{int}(F_{k})$ is dense} (see for example \cite[ Theorem 3.46, pg 93]{Hitchhiker}).
We will also use a result due to Avanissian
(\cite{Avanissian}, see also \cite{Lelong 1} ) according to which \textit{a separately subharmonic function is subharmonic if and only if it is locally bounded above.}

\vspace*{1cm}

 \section{Boundedness of a Family of Subharmonic Functions}
 
 \subsection{Singular set of a family of subharmonic functions} 
 
Let $\lbrace u_{\alpha}\rbrace_{\alpha\in J}$ be a family of subharmonic functions defined on an open set $\Omega$ of $\mathbb{R}^{N}$  and  set $u:=\sup_{\alpha\in J} u_{\alpha}$. Our first goal is to prove that the singular set of $u$, i.e., the set of points where $u^{*}$ is not subharmonic,  \textit{is a closed nowhere dense set whose components are all unbounded}. By definition $S$ is a closed set. The proof the $ S $ has empty interior is based on the Baire category theorem. However, since the functions $u_{\alpha}$ are not assumed to be continuous, we will encounter a technical problem. Following A. Sadullaev \cite{Sadullaev}, we will use properties of thin sets in classical potential theory to  circumvent this problem.

\begin{theorem}\label{th3.1.0}
Let $\Omega$ be a bounded open subset of $\mathbb{R}^{N}$  $(N\geq2)$ and  $(u_{\alpha})_{\alpha\in J}$   a family of subharmonic functions on a neighborhood of  $\overline{\Omega } $. If  $u(x):=\sup_{\alpha\in J}u_{\alpha}(x)$  is finite everywhere, then there exists a closed nowhere dense set $S\subset\Omega$ such that $u$ is locally bounded above at each point of $\Omega\setminus S$. Moreover, the set $S$ has no bounded components, and in particular, no isolated points. We will call $S$ the singular set of $u$.

\end{theorem}

\vspace*{10mm}

We need the following lemmas.
\begin{lemma}\label{lemma2.6.0}
Let  $ \omega $ be a bounded open subset of $ \mathbb{R}^{N} $ and $ v $ a subharmonic function defined on a neighborhood of $ \overline{\omega} $. Suppose that 
\begin{equation}\label{2.3.0} 
 v<\lambda 
\end{equation}
on $\omega$ for some constant $\lambda$. Then the set 
$$\Gamma=\left\lbrace \zeta\in \partial\omega:  v(\zeta)>\lambda \right\rbrace $$ is polar, and in particular,

$$v\leq \lambda$$
on $\mathrm{int}(\overline{\omega}).$
\end{lemma}

\begin{proof}
 It is sufficient  to prove that $\Gamma$ is polar.   We start by proving that $ \Gamma $  is negligible for $ \omega $, i.e., its harmonic measure is zero:
\begin{equation}
\mu_{x}^{\omega}(\Gamma)=0 \label{5555}
\end{equation} 
for all $ x\in \omega $.  To do so, it suffices to show that $ \omega$ is thin at each point  of $ \Gamma$, according to Theorem \ref{1000}-(i). Let $ \zeta\in \Gamma $. It follows from (\ref{2.3.0}) that
$$ \limsup_{\substack{x\rightarrow\zeta\\(x\in \omega)}}v(x)\leq \lambda $$
whereas $$ v(\zeta)>\lambda, $$ by definition of $ \Gamma $. Thus $ \omega$ is thin at $ \zeta $, according to  Theorem \ref{1002}. By Theorem \ref{1000} the set $ \Gamma $ is negligible  and (\ref{5555}) follows.

Next we proceed to prove that the set
$$\gamma:=\lbrace \zeta\in \partial \omega: v(\zeta)=\lambda\rbrace$$
is closed in $\partial \omega$. To do so it is sufficient to show that the restriction of $v$ to $\gamma$, that we steel write $v$, is continuous in the topological subspace $\partial \omega$. Since $v$ is already upper semi-continuous, we need to show that it is also lower semi-continuous in the mentioned topology. Let $\zeta\in \gamma$ and suppose that $v$ is not lower semi-continuous at $\zeta$. Then, there exits $\varepsilon>0$ such that for every open neighborhood $V$ of $\zeta$ there exists $x\in\left( \mathrm{int}(\overline{\omega})\cap\partial \omega\right) \cap V$ satisfying 
$$\lambda\leq v(x)\leq v(\zeta)-\varepsilon=\lambda-\varepsilon<\lambda,$$
which is absurd. Here, the first inequality is due to the fact that by (\ref{2.3.0}) the boundary of $\omega$ is included to the set 
\begin{equation}
\lbrace \xi:v(\xi)\geq\lambda\rbrace.\label{eq3.4}
\end{equation}
Thus the restriction of $v$ to $\gamma$ is continuous and $\gamma$ is closed in $\cap\partial \omega$.

Finally, let $ \chi $ be the characteristic function of $ \Gamma $ defined on $\partial\omega$. By (\ref{eq3.4}) we have $\partial\omega=\gamma\cup\Gamma$ and thus it follows form the last paragraph that $ \Gamma=\partial\omega\setminus\gamma $ is open in $ \partial\omega $; this implies   $\chi$ is continuous on $\Gamma$. Then, since  by (\ref{5555})  the harmonic extension of $ \chi $ satisfies $H_{\chi} \equiv0$, we obtain that all points of $\Gamma$ are irregular ( for the Dirichlet problem), and so is polar by Theorem\ref{1001}. For, if $\zeta\in\Gamma$ was regular, we would have 
$$0=\lim_{x\rightarrow\zeta}H_{\chi}(x)=\chi(\zeta)=1,$$ 
by continuity of $\chi$ on $\Gamma$. Impossible. 

\end{proof}

\vspace*{10mm}

\begin{lemma}[\textbf{The "maximum" principle}]\label{lemma3.2}
Let $\omega$ be a bounded domain of $\mathbb{R}^{N}$  $(N\geq2)$ and  $(u_{\alpha})_{\alpha\in J}$   a family of subharmonic functions on a neighborhood of   $\overline{\omega }$. Suppose there exist  a constant $M$ such that 
$$\sup_{\zeta\in\partial \omega}u(\zeta)\leq M,$$
where we set  $u:=\sup_{\alpha\in J}u_{\alpha}$. Then,  
\begin{equation}
\sup_{x\in\overline{\omega}} u(x)\leq M. \label{3.3.0}
\end{equation}

\end{lemma}
\begin{proof}
Let $\alpha\in J$. We have by hypothesis,
$$\limsup_{\substack{x\rightarrow\zeta\\(x\in\omega)}}  u_{\alpha}(x)\leq u_{\alpha}(\zeta)\leq u(\zeta)\leq M$$
for all $\zeta\in\partial\omega$. Thus by the maximum principle applied to $u_{\alpha}$, we have
$$u_{\alpha}(x)\leq M  $$
for all $x\in \omega$. This  implies (\ref{3.3.0}), since $\alpha$ is arbitrary.
\end{proof}
\vspace*{10mm}

\begin{proof}[Proof of Theorem \ref{th3.1.0} ]

   For $n=1,2,...$, let $E_{n} $ be the set of all $x\in \overline{\Omega}$ such that $u(x)<n$, and let $ E_{n,\alpha} $ be the set of all $x\in \overline{\Omega}$ such that $u_{\alpha}(x)<n$. Fix $ n $. It is clear that $ E_{n}\subset E_{n,\alpha} $  and then
\begin{equation}
\mathrm{int}(\overline{E_{n}})\subset \mathrm{int}(\overline{E_{n,\alpha}})\label{3.4}
 \end{equation}
for all $ \alpha\in J $.   
\vspace*{8mm}

 One can easily check that $\overline{\Omega}=\bigcup_{n=1}^{+\infty}\overline{E_{n}} $. Thus according to Baire  theorem, $G:=\bigcup_{n=1}^{+\infty}\text{int}(\overline{E_{n}})$ is dense in $\overline{\Omega}$.  Define $S:=\Omega\setminus G.$  To finish the proof of the first part of the theorem, it remains to show that $u$ is locally bounded above on  $\Omega\setminus S=G$.
  
   Let $\zeta\in G$.  There exists $n$ such that $\zeta\in \text{int}(\overline{E_{n}})$. Let $\alpha\in J$ be arbitrary.  By Lemma \ref{lemma2.6.0} applied to the subharmonic function $ u_{\alpha} $ on the open set $ E_{n,\alpha} $, we obtain   
 $$u_{\alpha}\leq n$$
  on $ \text{int}(\overline{E_{n,\alpha}}) $, and by (\ref{3.4}),  on $\text{int}(\overline{E_{n}})$.  Since $\alpha$ is arbitrary, 
  $$u\leq n,$$
on $ \text{int}(\overline{E_{n}}) $, which is an open neighborhood of $\zeta$. It follows that  $u$ is locally bounded above at $\zeta$, and since $\zeta$ is arbitrary, on $G$.  This prove the first part of the theorem.

\vspace*{8mm}

The second part of the theorem follows from  Lemma \ref{lemma3.2}. In fact, if $C$ were a bounded component of $S$, there would exist an open set $U$ containing $S$ without its boundary $\partial U$ meeting $S$. By the precited lemma, that would imply $U\cap S=\emptyset$. Absurd.
\end{proof}

\vspace*{10mm}
\remark
Recall that according to the (fundamental) convergence theorem, if $(u_{\alpha})$ is a family of locally  bounded-above subharmonic functions on a  domain $\Omega$ and $u:=\sup_{\alpha}u_{\alpha}$, then the upper regularization $u^{*}$ of $u$ is subharmonic and is equal to $u$ \textit{quasi-everywhere}  ( see for example \cite[pg 77]{Brelot} or \cite[pg 146]{Armitage Gardiner book} and the reference therein); here,
$$u^{*}(x):=\limsup_{y\rightarrow x}u(x).$$
If follows from Theorem 3.1  that the conclusion of the cited theorem  holds outside a nowhere dense closed set, if $u$ is finite everywhere. Also, according to Lemma 3.3, the same conclusion holds if $\Omega$ is bounded and $u$ is bounded above on $\partial\Omega$.

  \vspace*{10mm}

\subsection{The Converse Problem}

Given a set $S$ in $\mathbb{R}^{N}$, our goal is now to construct  a family of subharmonic functions whose supremum's singular set is $S$. By Theorem \ref{th3.1.0} it is necessary that  $S$  be closed nowhere dense with no bounded component. In our construction, we impose  an extra condition: that $\mathbb{R}^{N}\setminus S$ be connected.  We will follow closely  Wiegerinck's  example (\cite{Wiegerinck}, see also \cite{Wiegerinck Zeinstra}).

\begin{theorem} 
Let $S$ be a closed nowhere dense subset of $\mathbb{R}^{N}$ with no bounded component such that  $\mathbb{R}^{N}\setminus S$ is connected.  Then there exists   a family  $u_{\nu}$ of subharmonic functions in   $\mathbb{R}^{N}$ such that $u:=\sup_{\nu}u_{\nu}$ is finite  everywhere and locally unbounded at each point of $S$.
\end{theorem}

\begin{proof}
For $x\in\mathbb{R}^{N}$, let  $\text{d}(x,S):=\inf\lbrace \mid x-\zeta\mid: \zeta\in S\rbrace$. Following  Wiegerinck's example, we define for $n=1,2,...$,  
$$K_{n}:=\left\lbrace x\in \mathbb{R}^{N}:\mid x\mid\leq n,\hspace{1mm} \text{d}(x,S)\geq \frac{1}{n}\right\rbrace ,$$  
$$A_{n} :=\left\lbrace x\in \mathbb{R}^{N}:\mid x\mid\leq n,\hspace{1mm} \text{d}(x,S)= \frac{1}{n+1}\right\rbrace,$$
and
$$S_{n}:=S\cap\lbrace x:|x|\leq n\rbrace . $$
 Then we set
\begin{displaymath}
f_{n} = \left\{
\begin{array}{lr}
0 & \text{on a small neighborhood of  $K_{n}\cup S_{n}$}\\
n+1 &  \text{on a small neighborhood of $A_{n}$} 
\end{array}.
\right.
\end{displaymath}
This  is a harmonic function on a neighborhood of the compact $F_{n}:=K\cup A_{n}\cup S_{n}$. If follows from Theorem \ref{th2.7} that there exits a  harmonic function $H_{n} $ on $\mathbb{R}^{N}$  such that $\mid H_{n}\mid<\frac{1}{2}$ on $K_{n}\cup S$ and  $ \mid H_{n}\mid\geq n$ on $ A_{n}. $  Set $v_{n}=\max\lbrace \mid H_{n}\mid-1,0\rbrace$. Each $H_{n}$ being harmonic, $\mid H_{n}\mid$ is subharmonic everywhere and so is also $v_{n}$.  Clearly, the nonnegative function $v_{n}$ vanishes on $K_{n}\cup S_{n}$ and is bigger than $n$ on $A_{n}$.  Finally, we  define 
$$u_{\nu}(x)=\sum\limits_{n=1}^{\nu}v_{n}(x)$$ and $u(x)=\sup_{\nu}u_{\nu}(x)=\sum\limits_{n=1}^{+\infty}v_{n}(x).$ Notice that this sum is finite for all fixed $x\in\mathbb{R}^{N}$.

Now, take $x_{1}\in S$ and let us prove that $u$ is locally unbounded at $x_{1}$. Using the fact that $ S $ is closed with empty interior  one can easily construct a sequence $\lbrace\zeta_{n}\rbrace$  converging to $x_{1}$, as $n$ approaches infinity, and such that each $\zeta_{n}$ belongs to  $A_{n}$. 

To see this,  we first claim that for all $r>0$ there exists $N=N(r)$ such that $B(x_{1},r)$ encounters $ A_{n}$ for all $n\geq N$. Fix $r>0$ and take $y\in B(x_{1},r/4)\setminus S$ and $y'\in S$ such that $d(y,S)=d(y,y'):=\delta>0$. Then for $n\geq N>1/\delta$ there exists a point $\zeta_{n}$ in the segment $[y,y']$ joining $y$ to $y'$ such that $d(\zeta_{n},y')=1/n$. The choice of  $\zeta_{n}$ implies that it belongs to $A_{n}$ and also every point of $[y,y']$ belongs to $B(x_{1},r)$, and the claim follows.  Now using the claim one may easily see that the sequence $({\zeta_{n}}) $ converges to $x_{1}$, as $ n $ tends to infinity.

  Finally, we have $u(\zeta_{n})\geq\sum_{k=1}^{n}v_{k}(\zeta_{n})\geq n$, and thus $$\lim_{n\rightarrow+\infty}u(\zeta_{n})=+\infty.$$
\end{proof}

\vspace*{10mm}

\section{Application to Separately Subharmonic Functions}

Let $\Omega_{1}$ and $\Omega_{2}$ be two open sets in $\mathbb{R}^{p}$ and $\mathbb{R}^{q}$, respectively and $p,q\geq2$. A function $u(x,y)$ on $\Omega_{1}\times \Omega_{2}$ is called \textit{separately subharmonic} if for all fixed point $(a,b)\in \Omega_{1}\times \Omega_{2}$ the partial functions $x\mapsto u(x,b)$ and $y\mapsto u(a,y)$ are subharmonic on $\Omega_{1}$ and $ \Omega_{2}$, respectively. The following theorem generalizes a result due to Cegrell and Sadullaev  \cite[Theorem 6]{Cegrell Sadullaev} in which the authors assume moreover that $u(x,y)$ is  harmonic with respect to one variable.

\begin{theorem}
  Suppose that $u(x,y)$ is a separately subharmonic function on a neighborhood of $\overline{\Omega}_{1}\times\overline{\Omega}_{2}$. Then there exist two closed nowhere dense sets $S_{1}\subset \Omega_{1}$ and $S_{2}\subset \Omega_{2}$ such that $u(x,y)$ is subharmonic on $\Omega_{1}\times \Omega_{2}\setminus S_{1}\times S_{2}.$ Moreover, $S_{1}$ and $S_{2}$ have no bounded components, and in particular, no isolated points.
 \end{theorem}
 
 \begin{proof}
 Let $M(x):=\max_{\eta\in \overline{\Omega}_{2}}u(x,\eta)$. Since $u$ is upper semi-continuous with respect to the second variable, $M$ is a real-valued function. Thus we may apply Theorem 3.1, according to which there exists a closed nowhere dense  set $S_{1}\subset \Omega_{1}$ such that $M$ is locally bounded above at each point of $\Omega_{1}\setminus S_{1}$. Since $u\leq M$ on $\Omega_{1}\times \Omega_{2}$, if follows that $u(x,y)$ is locally bounded above at each point of $\Omega_{1}\setminus S_{1}\times \Omega_{2}$ and by Avanissian's theorem $u$ is subharmonic there. Same for $S_{2}$ by setting  $N(y):=\max_{\zeta\in \overline{\Omega}_{1}}u(\zeta,y)$. 
 \end{proof}

\begin{theorem}
Let $S_{1}\subset \mathbb{R}^{p}$ and $S_{2}\subset \mathbb{R}^{q}$ $(p,q\geq2)$ be two closed nowhere dense sets with no bounded component such that  $\mathbb{R}^{p}\setminus S_{1}$ and  $\mathbb{R}^{q}\setminus S_{2}$ are connected. Then there exists  a  function $u(x,y)$ separately subharmonic on $\mathbb{R}^{p}\times\mathbb{R}^{q}$, subharmonic on $ \mathbb{R}^{p}\times\mathbb{R}^{q}\setminus S_{1}\times S_{2}$, and locally unbounded at each point of $S_{1}\times S_{2}$.
\end{theorem}

\begin{proof}
As in the proof of Theorem 3.5,  define for $n=1,2,...,$ and $j\in\lbrace 1,2\rbrace$,  the  sets $K_{j,n}$, $A_{j,n}$ and $S_{j,n}$ and construct a subharmonic function $v_{j,n}(x)$ that vanishes on $K_{j,n}\cup S_{j,n}$ and is bigger that $n$ on $A_{j,n}$.  Then set $$u(x,y)=\sum\limits_{n=1}^{+\infty}v_{1,n}(x)v_{2,n}(y).$$
 If one of the variables $x$ or $y$ is fixed,  the above sum is finite and so the partial functions $u(.,y)$ and $u(x,.)$ are subharmonic. 

To see that $u$ is unbounded at a neighborhood of each point of $S_{1}\times S_{2}$, take a point $(x_{1},y_{1})$ in this set. As before, there exist two sequences $\lbrace \zeta_{n}\rbrace\subset A_{1,n}$ and $\lbrace \eta_{n}\rbrace\subset A_{2,n}$ converging to $ x_{1} $ and $ y_{1} $ respectively,  such that 
$$u(\zeta_{n},\eta_{n})\geq\sum\limits_{k=1}^{n}v_{1,k}(\zeta_{n})v_{2,k}(\eta_{n})\geq n^{2}.$$ Thus,
 $u(\zeta_{n},\eta_{n})\rightarrow+\infty$, as $n\rightarrow+\infty$. This proves that $u$ is locally unbounded at $(x_{1},y_{1})$ and so $S_{1}\times S_{2}$ is the singular set of $u$.
\end{proof}

\vspace*{10mm}
\remark The analog of Lemma \ref{lemma3.2}  for separately subharmonic functions shows that in order for such a function to be subharmonic, it is sufficient that it be bounded above on the distinguished boundary only. More precisely, 
\textit{Let $u(x,y)$ be a separately subharmonic function defined on a neighborhood of $\overline{\Omega_{1}\times\Omega_{2}}$, and assume  $ \Omega_{1} $ and $ \Omega_{2} $ are  connected. Suppose there exits a constant $M$ such that 
$$u\leq M$$
on $ \partial\Omega_{1}\times\partial\Omega_{2}$. Then, either $u<M$ or $u\equiv M$ on $\overline{\Omega_{1}\times\Omega_{2}}$, and in particular $u$ is subharmonic on $\Omega_{1}\times\Omega_{2}$.}

In fact, let $(x,y)\in \Omega_{1}\times\Omega_{2}.$ We have by the maximum principle applied to $u(x,.)$,
\begin{equation}\label{4.1}
u(x,y)\leq \max_{\eta\in\partial\Omega_{2}} u(x,\eta)=u(x,\eta_{1}),
\end{equation}
 for some $\eta_{1}\in\partial\Omega_{2}$. In the same way for $ \eta_{1} $ fixed, there exists $ \zeta_{1}\in\partial\Omega_{1} $ such that
$$u(x,\eta_{1})\leq \max_{\zeta\in\partial\Omega_{1}} u(\zeta,\eta)=u(\zeta_{1},\eta_{1}).$$
Thus  by hypothesis and in view of (\ref{4.1}),
$$u(x,y)\leq u(\zeta_{1},\eta_{1})\leq M,$$
 and so $u$ is bounded above on  $\Omega_{1}\times\Omega_{2}$.

Now, suppose $u(a,b)=M$ for some $(a,b)\in \Omega_{1}\times\Omega_{2}$ and take $(x,y)\in \Omega_{1}\times\Omega_{2}$. It is easy to see that
$$u(x,y)=M,$$
for $u(a,.)$ attains its maximum at $b\in\Omega_{2}$, and by the maximum principle, $u(a,.)\equiv M$; in particular, $u(a,y)=M$. Similarly, for $y$ fixed, $u(.,y)$ attains its maximum at $a\in\Omega_{1} $ and so is constant, in particular $u(x,y)=M$, as required.

 Finally, the last statement of the claim follows from Avanissian's theorem.

\end{document}